\documentclass[12pt,twoside]{article}
\usepackage[english]{babel}
\usepackage[latin1]{inputenc}
\usepackage{amsmath}
\usepackage{graphicx}
\usepackage{times,amssymb,amscd}
\newtheorem{thm}{Theorem}[section]

\newtheorem{defn}[thm]{Definition}
\newtheorem{rem}[thm]{Remark}
\numberwithin{equation}{section}

\newcommand{\bH}{\mathbf{H}}

\newcommand{\bL}{\mathbf{L}}

\newcommand{\bR}{\mathbf{R}}
\newcommand{\bS}{\mathbf{S}}

\newcommand{\br}{\mathbf{r}}

\newcommand{\bT}{\mathbf{T}}

\newcommand{\cP}{\mathcal{P}}
\newcommand{\cS}{\mathcal{S}}

\newcommand{\EUC}{\mathbf E^3}
\newcommand{\SPH}{\bS^3}
\newcommand{\HYP}{\bH^3}
\newcommand{\SXR}{\bS^2\!\times\!\bR}
\newcommand{\HXR}{\bH^2\!\times\!\bR}
\newcommand{\SLR}{\widetilde{\bS\bL_2\bR}}
\newcommand{\NIL}{\mathbf{Nil}}
\newcommand{\SOL}{\mathbf{Sol}}

\begin{document}
\pagestyle{myheadings}
\markboth{\centerline{Jen\H o Szirmai}}
{Desargues's and Pappus's hexagon theorems $\dots$}
\title
{Desargues's and Pappus's hexagon theorems on translation-type surfaces in Thurston geometries
}

\author{Jen\H o Szirmai \\
\normalsize Department of Algebra and Geometry, Institute of Mathematics,\\
\normalsize Budapest University of Technology and Economics, \\
\normalsize M\"uegyetem rkp. 3., H-1111 Budapest, Hungary \\
\normalsize szirmai@math.bme.hu
\date{\normalsize{\today}}}
\maketitle
\begin{abstract}
In \cite{Sz25} we generalized the famous Menelaus' and Ceva's theorems for translation triangles in each
non-constant curvature Thurston geometry. In this paper based on the described method and results, we prove  that the classical 
Desargues's and Pappus's hexagon theorems are true not only in classical geometries with constant curvature, but also in Thurston
geometries with non-constant curvature on the translation surfaces. In our work we use the unified projective models of Thurston geometries.

\end{abstract}
\newtheorem{Theorem}{Theorem}[section]
\newtheorem{corollary}[Theorem]{Corollary}
\newtheorem{lemma}[Theorem]{Lemma}
\newtheorem{exmple}[Theorem]{Example}
\newtheorem{definition}[Theorem]{Definition}
\newtheorem{rmrk}[Theorem]{Remark}
\newtheorem{proposition}[Theorem]{Proposition}
\newenvironment{remark}{\begin{rmrk}\normalfont}{\end{rmrk}}
\newenvironment{example}{\begin{exmple}\normalfont}{\end{exmple}}
\newenvironment{acknowledgement}{Acknowledgement}

\section{Introduction} \label{section1}
The role of non-Euclidean geometries in natural sciences is increasingly recognized, primarily in
physics. For all these reasons, three-dimensional maximal Riemannian geometries, also known as Thurston 
geometries, have come to the fore. The study of spaces with constant curvature
is very extensive, but the mapping of the internal structure of non-constant 
Thurston geometries still contains many open questions.
One of the most important area may be the non-Euclidean crystallography, since
it cannot be ruled out that under certain extreme conditions, materials
crystallize in ``different geometric structures´´.

The analysis of the structure of Thurston geometries 
was greatly facilitated by the description of the 
unified projective model of Thurston geometries 
by E. Moln\'ar in his article \cite{M97}. 
We use this model in the following.

Many classical theorems, which have been studied extensively in Euclidean geometry, 
can be extended to Thurston geometries with non-constant curvature. They offer new perspectives and applications in $\SXR, \HXR, \NIL, \SLR, \SOL$ spaces.

In the study of the internal structures of Thurston spaces, it is worth distinguishing 
two significantly different directions related to the translation distance or geodesic distance.

We can introduce in a natural way (see \cite{M97})
translation mappings any point to any other point. Consider a unit tangent
vector at the origin. Translations carry this vector to a tangent vector any
other point.
If a curve $t\rightarrow (x(t),y(t),z(t))$ has just the translated vector as its tangent vector at
each point, then the curve is called a {\it translation curve}. This assumption leads to
a system of first order differential equations. Thus, translation curves are simpler
than geodesics and differ from them in $\NIL$, $\SLR$ and $\SOL$ geometries.
In $\EUC$, $\SPH$, $\HYP$, $\SXR$ and $\HXR$ geometries, the translation and geodesic curves coincide
with each other. But in the $\NIL$, $\SLR$ and $\SOL$ geometries,
translation curves are in many ways more natural than geodesics. Therefore, we
distinguish two different distance functions: $d^g$ is the usual geodesic distance function,
and $d^t$ is the translation distance function. So we obtain
two types of curves, triangles, bisector surfaces (and two types of the corresponding Dirichlet-Voronoi cells) etc. by the two different distance functions,
but {\it{in the present paper we consider only the translation case}}. 

An important direction is the study of Apollonius surfaces, a special case of 
which is the bisector surfaces. These are of fundamental importance due to the determinations of Dirichlet-Voronoi cells.
The geodesic-like Apollonius surfaces are investigated in $\SXR$, $\HXR$ and $\NIL$ geometries in \cite{Sz22-2, Sz23-1} and in \cite{Cs-Sz25} 
the translation-like Apollonius surfaces in Thurston geometries.
Moreover, we studied the translation-like equidistant surfaces in $\SOL$ and $\NIL$ geometries in \cite{Sz19} and in \cite{VSz19}.
In \cite{PSSz10}, \cite{PSSz11-1}, \cite{PSSz11-2} we studied the geodesic-like bisector surfaces in $\SXR$ and $\HXR$ spaces.

We examined further interesting problems related to the Thurston geometries 
with non-constant curvature in papers \cite{CsSz16, Cs-Sz24, MSz12, MSzV17, Sz07, Sz10, Sz11, Sz12, 
Sz14-1, Sz14-2, Sz18, Sz20, Sz22-1, Sz23-2, Sz24}.

The next important question, which is also closely related to the
Apollonius surfaces, is the determination of the surface of a given geodesic or translational triangle. 
Defining this is an essential condition for stating elementary geometric concepts and theorems related to triangles.

In the works, \cite{Sz22-2}, \cite{Sz23-1} we proposed a possible definition of geodesic-like triangular surfaces in $\SXR$, $\HXR$ and $\NIL$ geometries, 
which was done with the help of Apollonius surfaces and returned traditional triangular surfaces in geometries with constant curvature.

\begin{rmrk}
However, defining the surface of a translation or geodesic triangle in Thurston geometries with non-constant curvature is not 
straightforward. The usual geodesic or translation-like triangle surface definition is not possible because 
the geodesic or translation curves starting from different vertices and ending at points of the corresponding opposite edges define different
surfaces in general, i.e. {\it geodesics or translation curves starting from different vertices and ending at points on the corresponding opposite side usually do not intersect.}
\end{rmrk}

In the paper \cite{Cs-Sz25}, we extended these questions to the translation triangles of non-constant curvature Thurston geometries. 
We provided a new possible definition of the surfaces of translation-like triangles and determined their equations.
In the following we will use these definitions of the surfaces of translation triangles in Thurston geometries.

Among the many classical results, Menelaus's theorem is particularly significant
in understanding external triangle properties and collinearity conditions.
While its formulation is well established in geometries of constant curvature, its behavior
in other geometries requires careful investigation. In Thurston geometries of constant curvature these theorems are well known and 
in \cite{Sz22-2} and \cite{Sz23-1} we generalized 
them to $\SXR$, $\HXR$ and $\NIL$ spaces. In \cite{Sz25} we continued to investigate this issue, 
generalized and proved the theorems of Menelaus' and Ceva's for {\it translation triangles} in $\NIL$, $\SLR$ and $\SOL$ spaces. 

However, a systematic
extension of the consequences of the Ceva's and Menelaus' theorems and the related projective geometry theorems to
alternative metric spaces has not been thoroughly examined.

{\it This paper aims to fill this gap by analyzing Pappus's theorem and Desargues's theorem in Thurston geometries with non-constant curvatures.
By establishing these theorems in different geometric structures,
this study provides new insights into the relationships between Euclidean
and non-Euclidean geometries.}

In the structure of the article, we follow the principle of giving a complete overview of the proof in the case of $\NIL$ 
geometry, and in the case of other geometries, since proofs with similar structures occur, we only cover the space-specific differences.

\section{Pappus's and Desargues's theorems for translation triangles in $\NIL$ geometry}
\subsection{On the projective model of $\NIL$ geometry}
$\NIL$ geometry can be derived from the famous real matrix group $\mathbf{L(\mathbf{R})}$ discovered by Werner Heisenberg. The left (row-column) 
multiplication of Heisenberg matrices
     \begin{equation}
     \begin{gathered}
     \begin{pmatrix}
         1&x&z \\
         0&1&y \\
         0&0&1 \\
       \end{pmatrix}
       \begin{pmatrix}
         1&a&c \\
         0&1&b \\
         0&0&1 \\
       \end{pmatrix}
       =\begin{pmatrix}
         1&a+x&c+xb+z \\
         0&1&b+y \\
         0&0&1 \\
       \end{pmatrix}
      \end{gathered} \label{2.7}
     \end{equation}
defines "translations" $\mathbf{L}(\mathbf{R})= \{(x,y,z): x,~y,~z\in \mathbf{R} \}$ 
on the points of $\NIL= \{(a,b,c):a,~b,~c \in \mathbf{R}\}$. 
These translations are not commutative in general. The matrices $\mathbf{K}(z) \vartriangleleft \mathbf{L}$ of the form
     \begin{equation}
     \begin{gathered}
       \mathbf{K}(z) \ni
       \begin{pmatrix}
         1&0&z \\
         0&1&0 \\
         0&0&1 \\
       \end{pmatrix}
       \mapsto (0,0,z)  
      \end{gathered}\label{2.8}
     \end{equation} 
constitute the one parametric centre, i.e. each of its elements commutes with all elements of $\mathbf{L}$. 
The elements of $\mathbf{K}$ are called {\it fibre translations}. $\NIL$ geometry of the Heisenberg group can be projectively 
(affinely) interpreted by "right translations" 
on points as the matrix formula 
     \begin{equation}
     \begin{gathered}
       (1;a,b,c) \to (1;a,b,c)
       \begin{pmatrix}
         1&x&y&z \\
         0&1&0&0 \\
         0&0&1&x \\
         0&0&0&1 \\
       \end{pmatrix}
       =(1;x+a,y+b,z+bx+c) 
      \end{gathered} \label{2.9}
     \end{equation} 
shows, according to (\ref{2.7}). Here we consider $\mathbf{L}$ as projective collineation group with right actions in homogeneous coordinates.
We will use the usual projective model of $\NIL$ (see \cite{M97}).

The translation group $\mathbf{L}$ defined by formula (\ref{2.9}) can be extended to a larger group $\mathbf{G}$ of collineations,
preserving the fibres, that will be equivalent to the (orientation preserving) isometry group of $\NIL$. 

In \cite{MSz} we has shown that 
a rotation through angle $\omega$
about the $z$-axis at the origin, as isometry of $\NIL$, keeping invariant the Riemann
metric everywhere, will be a quadratic mapping in $x,y$ to $z$-image $\overline{z}$ as follows:
     \begin{equation}
     \begin{gathered}
       \mathcal{M}=\br(O,\omega):(1;x,y,z) \to (1;\overline{x},\overline{y},\overline{z}); \\ 
       \overline{x}=x\cos{\omega}-y\sin{\omega}, \ \ \overline{y}=x\sin{\omega}+y\cos{\omega}, \\
       \overline{z}=z-\frac{1}{2}xy+\frac{1}{4}(x^2-y^2)\sin{2\omega}+\frac{1}{2}xy\cos{2\omega}.
      \end{gathered} \label{2.10}
     \end{equation}
This rotation formula $\mathcal{M}$, however, is conjugate by the quadratic mapping $\alpha$ to the linear rotation $\Omega$ as follows
     \begin{equation}
     \begin{gathered}
       \alpha^{-1}: \ \ (1;x,y,z) \stackrel{\alpha^{-1}}{\longrightarrow} (1; x',y',z')=(1;x,y,z-\frac{1}{2}xy) \ \ \text{to} \\
       \Omega: \ \ (1;x',y',z') \stackrel{\Omega}{\longrightarrow} (1;x",y",z")=(1;x',y',z')
       \begin{pmatrix}
         1&0&0&0 \\
         0&\cos{\omega}&\sin{\omega}&0 \\
         0&-\sin{\omega}&\cos{\omega}&0 \\
         0&0&0&1 \\
       \end{pmatrix}, \\
       \text{with} \ \ \alpha: (1;x",y",z") \stackrel{\alpha}{\longrightarrow}  (1; \overline{x}, \overline{y},\overline{z})=(1; x",y",z"+\frac{1}{2}x"y").
      \end{gathered} \label{2.11}
     \end{equation}
This quadratic conjugacy modifies the $\NIL$ translations in (2.3), as well. 
\subsubsection{Translation curves}
We consider a $\NIL$ curve $(1,x(t), y(t), z(t) )$ with a given starting tangent vector at the origin $O=E_0=(1,0,0,0)$
\begin{equation}
   \begin{gathered}
      u=\dot{x}(0),\ v=\dot{y}(0), \ w=\dot{z}(0).
       \end{gathered} \label{2.12}
     \end{equation}
For a translation curve let its tangent vector at the point $(1,x(t), y(t), z(t) )$ be defined by the matrix (\ref{2.9}) 
with the following equation:
\begin{equation}
     \begin{gathered}
     (0,u,v,w)
     \begin{pmatrix}
         1&x(t)&y(t)&z(t) \\
         0&1&0&0 \\
         0&0&1&x(t) \\
         0&0&0&1 \\
       \end{pmatrix}
       =(0,\dot{x}(t),\dot{y}(t),\dot{z}(t)).
       \end{gathered} \label{2.13}
     \end{equation}
Thus, the {\it translation curves} in $\NIL$ geometry (see  \cite{Cs-Sz23}, \cite{MSzV}) are defined by the above first order differential equation system 
$\dot{x}(t)=u, \ \dot{y}(t)=v,  \ \dot{z}(t)=v \cdot x(t)+w,$ whose solution is the following: 
\begin{equation}
   \begin{gathered}
       x(t)=u t, \ y(t)=v t,  \ z(t)=\frac{1}{2}uvt^2+wt.
       \end{gathered} \label{2.14}
\end{equation}
We assume that the starting point of a translation curve is the origin, because we can transform a curve into an 
arbitrary starting point by translation (2.3), moreover, unit initial velocity translation can be assumed by "geographic" parameters $\phi$ and $\theta$:
\begin{equation}
\begin{gathered}
        x(0)=y(0)=z(0)=0; \\ \ u=\dot{x}(0)=\cos{\theta} \cos{\phi}, \ \ v=\dot{y}(0)=\cos{\theta} \sin{\phi}, \ \ w=\dot{z}(0)=\sin{\theta}; \\ 
        - \pi \leq \phi \leq \pi, \ -\frac{\pi}{2} \leq \theta \leq \frac{\pi}{2}. \label{2.15}
\end{gathered}
\end{equation}
\subsection{The surfaces of translation triangles and generalizations of Menelaus' theorem}
We consider $3$ points $A_0$, $A_1$, $A_2$ in the projective model of the space $X$ ( $X \in \{\EUC,\SPH,\HYP,\SXR,\HXR,\NIL,\SLR,\SOL\}$).
The {\it translation segments} $a_k$ connecting the points $A_i$ and $A_j$
$(i<j,~i,j,k \in \{0,1,2\}, k \ne i,j$) are called sides of the {\it translation triangle} with vertices $A_0$, $A_1$, $A_2$. 

We introduced a new definition of the surface $\mathcal{S}_{A_0A_1A_2}$ of the translation triangle in \cite{Cs-Sz25}: 
{\it The $S^{X,t}_{A_0A_1A_2}$ translation-like triangular surface of translation triangle $A_0A_1A_2$ in the Thurston geometry $X$ is  
the set of all $P$ points of $X$ from which the tangents in $P$ of the translation curves drawn to vertices $A_0$ $A_1$ and $A_2$ are coplanar in Euclidean sense.}
We also determined the equations of these surfaces and visualized them and we will use this surface definition in the following.
An important question is how to define the curves that pass through the surface of a translation triangle. In the $\EUC,\SPH,\HYP,\NIL,\SLR,\SOL$ geometries, 
the surface defined above contains the translation curve connecting any two of its points, but not in the $\SXR,\HXR$ geometries. 
Therefore, we discussed its definition separately in these geometries (see \cite{Sz22-2}, \cite{Sz25}).
The Ceva's and Menelaus' theorems were already known in geometries of constant curvature $\EUC, \SPH, \HYP$. 
We generalized them to geodesic triangles in $\SXR, \HXR$ and $\NIL$ geometries (see \cite{Sz22-2}, \cite{Sz23-1} ) and to all 
Thurston geometries for translational triangles in \cite{Sz25}. Now, we recall the result for $\NIL$ geometry, related to the Menelaus' theorem.

We consider a {\it translation triangle}
$A_0A_1A_2$ in the projective model of the $\NIL$ space 
Without limiting generality, we can assume that $A_0=(1,0,0,0)$.
The translation curves that contain the sides $A_0A_1$ and $A_0A_2$ of the given triangle can be characterized directly by the corresponding parameters $\theta$ 
and $\phi$ (see (2.8) and (2.9)).

The translation curve including the translation-like side segment $A_1A_2$ is also determined by one of its endpoints and its parameters 
but in order to determine the corresponding parameters of this
translation curve we use the translation  $\bT^{\NIL}(A_1)$, as elements of the isometry group of the geometry $\NIL$, that
maps $A_1=(1,x_1,y_1,z_1)$ onto $A_0=(1,0,0,0)$. 

From the (2.8), (2.9) equations of the $\NIL$ translation curve we directly obtain the following
\begin{lemma}
\begin{enumerate}
\item Let $P_1$ and $P_2$ be arbitrary points and $g_{P_1,P_2}^\NIL$ is the corresponding translation curve in the considered model of $\NIL$ geometry. 
The projection of the translation curve segment $g_{P_1,P_2}^\NIL$ onto the $[x,y]$ coordinate plane 
in the direction of the fibers will be a segment in the Euclidean sense if the points of the translation curve do not lie on a fiber
otherwise, its projection is a point.

\item If points $P_1$ and $P_2$ lie in the $[x,z]$ or $[y,z]$ coordinate plane, then the corresponding translation curve
$g_{P_1,P_2}^{\NIL}$ is a straight line segment $P_1P_2$ in the Euclidean sense.
\end{enumerate}~ ~ $\square$
\end{lemma}
\begin{rmrk}
In other words, the points of the translation curve $g_{P_1,P_2}^\NIL$ lie exactly in a plane in Euclidean sense orthogonal to the $[x,y]$ coordinate plane
if the points of the translation curve do not lie on a fiber.
\end{rmrk}
We extend the definition of the simple ratio to the $\NIL$ space for translation curves.
\begin{definition}[\cite{Sz25}]
If $A$, $B$ and $P$ are distinct points on a translation curve in the $\NIL$ space, then
their simple ratio is
$$s_g^\NIL(A,P,B) =  d^{\NIL,t}(A,P)/d^{\NIL,t}(P,B),$$ if $P$ is between $A$ and $B$, and
$$s_g^\NIL(A,P,B) = -d^{\NIL,t}(A,P)/d^{\NIL,t}(P,B),$$.
\end{definition}
\begin{lemma}[\cite{Sz25}]
\begin{enumerate}
\item Let $g_{A,B}$ be an arbitrary non-fibrum-like translation curve, where without loss of generality we can assume 
that $A=(1,0,0,0)$ coincides with the center of the model. 
Furthermore, let $P \in g^\NIL_{A,B}$ and let $A$, $B$ be the projected images 
by fibers into the $[x,y]$ coordinate plane, $A^p$, $B^p$, and $P^p$ then
$$
s_g^\NIL(A,P,B)=s_g^\NIL(A^p,P^p,B^p)=s_g^{\EUC}(A^p,P^p,B^p),
$$ 
\item If points $A$, $B$ and $P$ lie on a fibrum, then
$$
s_g^\NIL(A,P,B)=s_g^{\EUC}(A,P,B)
$$
\end{enumerate}
\end{lemma} \quad \quad $\square$

We introduced the following notations:

1. If the surface of a translation-like triangle is a plane in Euclidean sense
then it is called fibre type triangle.

2. In the other cases the triangle is in general type (see \cite{Cs-Sz25}).
\begin{definition}[\cite{Sz25}]
Let $\mathcal{S}^{\NIL,t}_{A_0A_1A_2}$ be the surface of the translation triangle $A_0A_1A_2$ and $P_1$, $P_2 \in \mathcal{S}^{\NIL,t}_{A_0A_1A_2}$ two given point.
the connecting curve $P_1P_2 \subset \mathcal{S}^{\NIL,t}_{A_0A_1A_2}$ is a translation curve in $\NIL$ space.
\end{definition}
\begin{Theorem}[Menelaus's theorem for translation triangles in $\NIL$ space \cite{Sz25}]
If $l$ is a line not through any vertex of a translation triangle
$A_0A_1A_2$ lying in a surface $\mathcal{S}^{\NIL,t}_{A_0A_1A_2}$
such that
$l$ meets the translation curves $g_{A_1A_2}^X$ in $Q$, $g_{A_0A_2}^X$ in $R$,
and $g_{A_0A_1}^\NIL$ in $P$,
then $$s^\NIL_g(A_0,P,A_1)s^\NIL_g(A_1,Q,A_22)s^\NIL_g(A_2,R,A_0) = -1.$$
\end{Theorem}~ ~ $\square$
\begin{rmrk}
It is easy to see that the ``reversals'' of the above theorems are also true.
\end{rmrk}
\subsection{Desargues's theorem in $\NIL$ geometry} 
The theorem of Desargues are the fundamental
building blocks in the axiomatic development of incidence and projective
geometry. Desargues's theorem is a fundamental theorem in projective geometry demonstrates that corresponding sides and vertices
of two triangles are projectively connected. An important aspect of this theorem
is that real projective geometry extends beyond plane geometry to encompass
spatial geometry. Desargues's theorem is essential in explaining the principles
of perspective and is widely used in applications of projective geometry.  
In dimension $2$ the planes for which it holds are called {\it Desarguesian planes}. 

Desargues's theorem, in geometry, mathematical statement discovered by the French mathematician Girard Desargues in 1639 
that motivated the development, in the first quarter of the 19th century, of projective geometry.

First, let's consider the following original interpretation of Desargues' theorem in the $\NIL$ space.
\begin{figure}[ht]
\centering
\includegraphics[width=12cm]{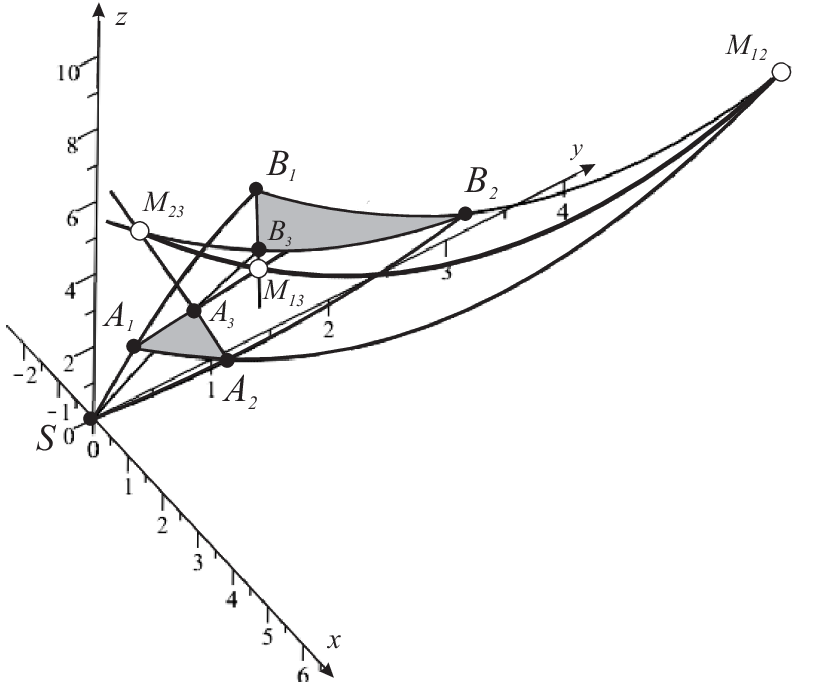}
\caption{Desargues's theorem on translation-type surface in $\NIL$ geometry.}
\label{}
\end{figure}
\begin{defn}
Consider a translation triangle $A_1A_2A_3$ in the space $\NIL$ whose translation surface defined in subsection 2.2 (see \cite{Cs-Sz25}), is denoted 
by $\cS_{A_1A_2A_3}^{\NIL,t}$. Moreover, let the translation triangle $B_1B_2B_3$ be located in the same surface. We assume for the above two triangles that 
the translation curves 
$g_{A_1B_1}^\NIL$, $g_{A_2B_2}^\NIL$, $g_{C_1C_2}^\NIL$ intersect at the proper point $S \subset \cS_{A_1A_2A_3}^{\NIL,t}$. 
Then we say that the triangles $A_1A_2A_3$ and $B_1B_2B_3$ are {\it perspective with respect to the point, (perspective centrally)} in the 
$\NIL$ sense. Moreover, if for the initial triangles it is true that ${M_{12}}=g_{A_1A_2}^\NIL \cap g_{B_1B_2}^\NIL$, $M_{13}=g_{A_1A_3}^\NIL \cap 
g_{B_1B_3}^\NIL$, $M_{23}=g_{A_2A_3}^\NIL \cap g_{B_2B_3}^\NIL$
proper points exist and fit on a translation curve, the triangles are {\it perspective with respect to a line, (perspective axially)} in the $\NIL$ space.
\end{defn}
\begin{Theorem}[Desargues's theorem in $\NIL$ space] If two triangles $A_1A_2A_3$ and $B_1B_2B_3$ are perspective centrally with respect to the point $S$, 
then the points $M_{12}$, $M_{13}$ and $M_{23}$ of intersection of the two triangles' corresponding
sides (translation curves) - if those intersection points exist - all lie on one translation curve.
\end{Theorem}

{\bf{Proof}}

Consider translation triangle $A_1A_2A_3 \subset \NIL$ (see Fig.~1) and a point $S \subset \cS_{A_1A_2A_3}^{\NIL,t}$ different from its vertices.
Let $B_1 \ne A_1$ lies on $g_{SA_1}^\NIL$,  $B_2 \ne A_2$ lies on 
$g_{SA_2}^\NIL$, $B_3 \ne A_3$ lies on $g_{SA_3}^\NIL$.

\begin{enumerate}
\item Let us apply the $\NIL$ version of Menelaus' theorem (Theorem 2.6) 
to the translation triangle $SA_3A_1$ and the collinear points $B_1, B_3, M_{13}$:
\begin{equation}
s^\NIL_g(S,B_1,A_1)s^\NIL_g(A_1,M_{13},A_3)s^\NIL_g(S,B_3,A_3) = -1. \tag{2.10}
\end{equation}
 \item Similarly, applying the $\NIL$ version of Menelaus' theorem for translation
triangle $SA_2A_3$ and the collinear points $B_2, B_3, M_{23}$:

\begin{equation}
s^\NIL_g(S,B_3,A_3)s^\NIL_g(A_3,M_{23},A_2)s^\NIL_g(S,B_2,A_2) = -1. \tag{2.11}
\end{equation}
\item Similarly to the above cases, by the $\NIL$ version of Menelaus' theorem for translation triangle $SA_1A_2$ 
and the collinear points $B_1, B_2, M_{12}$:
\begin{equation}
s^\NIL_g(S,B_1,A_1)s^\NIL_g(A_1,M_{12},A_2)s^\NIL_g(S,B_2,A_2) = -1. \tag{2.12}
\end{equation}
\item By multiplying all of the left sides together and all of the right sides together, we
find
\begin{equation}
\begin{gathered}
s^\NIL_g(S,B_1,A_1)s^\NIL_g(A_1,M_{13},A_3)s^\NIL_g(S,B_3,A_3) \cdot \\
\cdot s^\NIL_g(S,B_3,A_3)s^\NIL_g(A_3,M_{23},A_2)s^\NIL_g(S,B_2,A_2) \cdot \\
\cdot s^\NIL_g(S,B_1,A_1)s^\NIL_g(A_1,M_{12},A_2)s^\NIL_g(S,B_2,A_2) = \\
\frac{A_{1}M_{13}}{M_{13}A_{3}} \cdot \frac{A_{3}M_{23}}{M_{23}A_{2}} 
\cdot \frac{A_{2}M_{12}}{M_{12}A_{1}} = -1= \\
=s^\NIL_g(A_1,M_{13},A_3)s^\NIL_g(A_3,M_{23},A_2)s^\NIL_g(A_2,M_{12},A_1)=-1. \tag{2.13}
\end{gathered}
\end{equation}

\end{enumerate}
Applying the reversal of the $\NIL$ version of the Menelaus' theorem to 
the translation triangle $A_1A_2A_3$, we obtain that the points
$M_{12}, M_{13}, M_{23}$ lie on a translation curve. 
collinear. \qquad \qquad $\square$
\subsection{Pappus's hexagon theorem in $\NIL$ geometry} 
Pappus's hexagon theorem (attributed to Pappus of Alexandria), similarly to the 
Desargues's theorem, it also plays a fundamental 
role in the axiomatic structure of projective geometry. Projective planes in which 
the theorem is valid are called {\it Pappian planes}. 

Consider a translation plane $\mathcal{S}^{\NIL,t}_{g_a,g_b}$ given two translation curves $g_a^\NIL$ and $g_b^\NIL$ where points $A_1, A_2, A_3$ 
are different points on the translation curve $g_a^\NIL$ and $B_1, B_2, B_3$ 
are different points on the translation curve $g_b^\NIL$ (see Fig.~3).

\begin{Theorem}[Pappus's hexagon theorem in $\NIL$ space]  

Let $A_1, A_2, A_3 \subset \mathcal{S}^{\NIL,t}_{g_a,g_b}$  
be three points on a translation curve $g_a^\NIL$ and let $B_1, B_2, B_3 
\subset \mathcal{S}^{\NIL,t}_{g_a,g_b}$ be three points on another translation 
curve $g_b^\NIL$. 
If the translation curves $g_{A_1B_1}^\NIL$, 
$g_{A_1B_1}^\NIL$, $g_{A_2B_2}^\NIL$, $g_{A_3B_3}^\NIL$ intersect the translation 
curves $g_{A_3B_2}^\NIL$,  $g_{A_1B_3}^\NIL$, $g_{A_2B_1}^\NIL$, respectively then 
then the three points of intersection $M_1=g_{A_1B_1}^\NIL \cap g_{A_3B_2}^\NIL$,
$M_2=g_{A_2B_2}^\NIL \cap g_{A_1B_3}^\NIL$, 
$M_3=g_{A_3B_3}^\NIL \cap g_{A_2B_1}^\NIL$ - if those 
intersection points exist - all lie on one translation curve (see Fig.~3).
\end{Theorem}

{\bf{Proof}}

We introduce the following notation: 
$C_{12}=g_{A_1B_1}^\NIL \cap g_{A_2B_2}^\NIL$, 
$C_{13}=g_{A_1B_1}^\NIL \cap g_{A_3B_3}^\NIL$ and 
$C_{23}=g_{A_2B_2}^\NIL \cap g_{A_3B_3}^\NIL$.

Consider the translation triangle $C_{12}C_{13}C_{23} \subset \mathcal{S}^{\NIL,t}_{g_a,g_b}$.

\begin{enumerate}
\item Let us apply the $\NIL$ version of Menelaus' theorem (Theorem 2.6) 
to the translation triangle $C_{12}C_{13}C_{23}$ and the collinear points $A_3, B_2, M_{1}$:
\begin{equation}
s^\NIL_g(C_{12},B_2,C_{23})s^\NIL_g(C_{23},A_{3},C_{13})s^\NIL_g(C_{13},M_1,C_{12}) =
-1. \tag{2.14}
\end{equation}
 \item Similarly, applying the $\NIL$ version of Menelaus' theorem for translation
triangle $C_{12}C_{13}C_{23}$ and the collinear points $A_1, B_3, M_{2}$:
\begin{equation}
s^\NIL_g(C_{12},M_2,C_{23})s^\NIL_g(C_{23},B_{3},C_{13})s^\NIL_g(C_{13},A_1,C_{12}) 
= -1.\tag{2.15}
\end{equation}
\item Similarly to the above cases, by the $\NIL$ version of Menelaus' 
theorem for translation triangle $C_{12}C_{13}C_{23}$ 
and the collinear points $B_1, A_2, M_{3}$:
\begin{equation}
s^\NIL_g(C_{12},A_2,C_{23})s^\NIL_g(C_{23},M_{3},C_{13})s^\NIL_g(C_{13},B_1,C_{12}) = 
-1. \tag{2.16}
\end{equation}
\item Similarly to the above cases, by the $\NIL$ version of Menelaus' 
theorem for translation triangle $C_{12}C_{13}C_{23}$ 
and the collinear points $A_1, A_2, A_{3}$:
\begin{equation}
s^\NIL_g(C_{12},A_2,C_{23})s^\NIL_g(C_{23},A_{3},C_{13})s^\NIL_g(C_{13},A_1,C_{12}) = 
-1. \tag{2.17}
\end{equation}
\item Similarly to the above cases, by the $\NIL$ version of Menelaus' 
theorem for translation triangle $C_{12}C_{13}C_{23}$ 
and the collinear points $B_1, B_2, B_{3}$:
\begin{equation}
s^\NIL_g(C_{12},B_2,C_{23})s^\NIL_g(C_{23},B_{3},C_{13})s^\NIL_g(C_{13},B_1,C_{12}) = -1. \tag{2.18}
\end{equation}
\item By multiplying all of the left sides together and all of the right sides together, we
find
\begin{equation}
\begin{gathered}
\frac{C_{12}M_2}{M_2C_{23}} \cdot \frac{C_{23}M_3}{M_3C_{13}} 
\cdot \frac{C_{13}M_1}{M_1C_{12}} = -1= \\
=s^\NIL_g(C_{12},M_2,C_{23})s^\NIL_g(C_{23},M_{3},C_{13})s^\NIL_g(C_{13},M_1,C_{12}) = -1.
\tag{2.19}
\end{gathered}
\end{equation}

\end{enumerate}
Applying the reversal of the $\NIL$ version of the Menelaus' theorem to 
the translation triangle $A_1A_2A_3$, we obtain that the points $M_{12}, M_{13}, M_{23}$ lie in a translation curve 
``collinear" in the $\NIL$ translation sense. \qquad \qquad $\square$
\begin{figure}[ht]
\centering
\includegraphics[width=10cm]{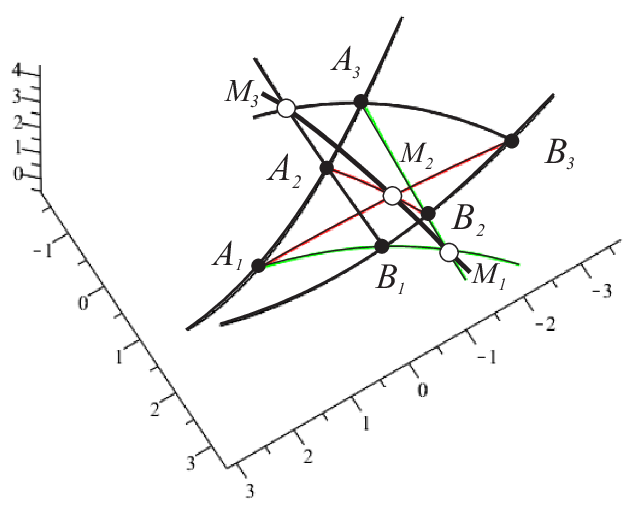}
\caption{Pappus's hexagon theorem on translation-type surface in $\NIL$ geometry.}
\label{}
\end{figure}
\begin{figure}[ht]
\centering
\includegraphics[width=8cm]{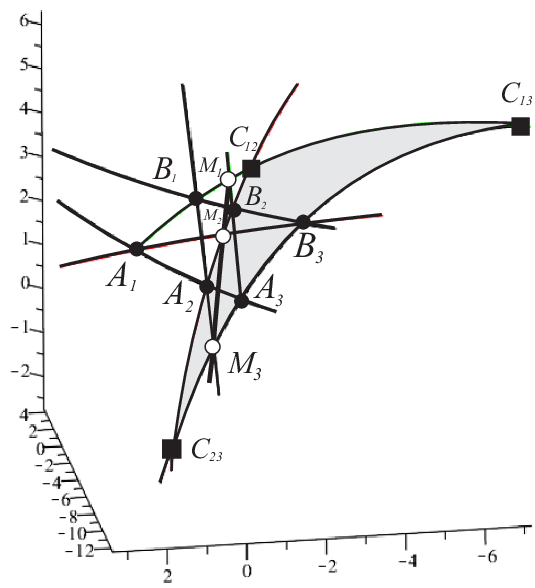}
\caption{The proof of the Pappus's hexagon theorem using Menelaus' theorem on $\NIL$ translation-type surface.}
\label{}
\end{figure}
\subsection{The extensions of the Desargues's and Pappus's hexagon theorems}
It is well known that the Euclidean plane can be embedded in, or extended to the real projective plane
as follows. We adjoin an ``infinite point" to each
of a given parallel class of lines. Different parallel classes have different infinite
points. We decree that the newly created infnite points lie on an ``infnite
line". The projective plane now has the property that
two distinct points lie on a unique line and two distinct lines meet in a
unique point.

A similar extension can be introduced for a given $\NIL$ translation surface introduced by us in the work \cite{Cs-Sz25}. Let there be a translation surface 
$\cS_{P_1P_2P_3}^{\NIL,t}$ with three points $P_1,P_2,P_3$ that are not located on a translation curve and do not lie in the Euclidean coordinate 
planes $[x,z]$ or $[y,z]$ or in a Euclidean plane parallel to them.
\begin{rem}
If the points are in the $[x,z]$, $[y,z]$ coordinate planes or in Euclidean 
planes parallel to them, then we get the Euclidean-type case, thus the Desargues's and Pappus's hexagon theorems are automatically satisfied
(see Lemmas 2.1, 2.4 and Theorem 2.6).
\end{rem}
On this surface, two {\it translation curves are called parallel} if they do not have a real common intersection point.
This is the same as the Euclidean planes that contain the translation curves and are parallel to the $z$-axis of the model are parallel in the Euclidean sense. 
All of the members of a set of mutually parallel translation curves have the same single point at infinity in common, 
and for each different orientation (belonging to different ``vertical Euclidean planes)
that contain parallel translation curves, there is a distinct point at infinity. All (and only)
the points at infinity arising from translation curves on a given translation surface lie on a given ``translation curve at
infinity". Therefore, this generalized projective translation surface $\cS\cP_{P_1P_2P_3}^{\NIL,t}$ consists of the given translation surface $\cS_{P_1P_2P_3}^{\NIL,t}$
plus a translation curve $g_{\cS_{P_1P_2P_3}^{\NIL,t}}^i $ at infinity. Considering this extension and the 
Remark 2.2, the Desargues's and Pappus's hexagon theorems can be easily extended by taking into account the Euclidean analogy.

\begin{Theorem}[Projective Desargues's theorem in $\NIL$ space] 

If two triangles $A_1A_2A_3$ and $B_1B_2B_3$ are perspective centrally (see Definition 2.1) with respect to the point $S$, 
then the points $M_{12}$, $M_{13}$ and $M_{23}$ of intersection of the two triangles' corresponding
sides (translation curves) all lie on one translation curve (perspective axially, see Definition 2.1) (see Fig.~1). \qquad \qquad $\square$
\end{Theorem} 
\begin{Theorem}[Projective Pappus's hexagon theorem in $\NIL$ space] 

Let $A_1$, $A_2$, $A_3 \subset \mathcal{S}^{\NIL,t}_{g_a,g_b}$  
be three points on a translation curve $g_a^\NIL$ and let $B_1, B_2, B_3 
\subset \mathcal{S}^{\NIL,t}_{g_a,g_b}$ be three points on another translation 
curve $g_b^\NIL$. 
The translation curves $g_{A_1B_1}^\NIL$, 
$g_{A_1B_1}^\NIL$, $g_{A_2B_2}^\NIL$, $g_{A_3B_3}^\NIL$ intersect the translation 
curves $g_{A_3B_2}^\NIL$,  $g_{A_1B_3}^\NIL$, $g_{A_2B_1}^\NIL$, respectively then 
then the three points of intersection $M_1=g_{A_1B_1}^\NIL \cap g_{A_3B_2}^\NIL$,
$M_2=g_{A_2B_2}^\NIL \cap g_{A_1B_3}^\NIL$, 
$M_3=g_{A_3B_3}^\NIL \cap g_{A_2B_1}^\NIL$ all lie on one translation curve (see Fig.~2,3). \qquad \qquad $\square$
\end{Theorem} 
\section{Extension of the theorems of Pappus and Desargues to additional Thurston geometries}
It is well known that the use of the projective models of non-euclidean geometries of constant curvature has several advantages. 
Many projective theorems could have multiple interpretations as theorems in elliptic or hyperbolic space geometry. 
It is also obvious that the theorems of hyperbolic space and the plane can be derived from many 
projective theorems, but several difficulties arise here, 
let us just think that in hyperbolic geometry, lines can be ultraparallel, i.e. they intersect outside the model.
Most of all this can be eliminated by introducing generalized concepts, such as the concept of a generalized triangle 
by introducing the polar line, but it is clear that these concepts require further investigation. 

Similar problems arise with some of the Thurston geometries with non-constant curvature. This is clearly shown by the 
problems that arise in the generalization of the theorems of Ceva and Menelaus (see \cite{Sz25}), which are related to the 
structure of projective models of corresponding Thurston geometries (see \cite{M97}, \cite{MSz}). Because of all this, we can 
also formulate the classical form of the Desargues's and Pappus's hexagon theorems. 
Their extension requires further investigation.

In the In Thurston geometries of constant curvature $\EUC, \SPH, \HYP$ the Desargues's and Pappus's hexagon theorems are well 
known and have been discussed, see for example \cite{AB13}, \cite{CG}.
In the previous section we saw the Desargues's and Pappus's hexagon theorems in $\NIL$ geometry. 

Now, we consider the same theorems in the further non-constant curvature Thurston geometries. We will not detail the proofs of these theorems, since 
they can be done similarly to $\NIL$ geometry (see Section 2), because Menelaus' theorem holds in these geometries as well (see \cite{Sz25}).
\begin{defn}
Let $X \in \{\SXR, \HXR, \SOL, \SLR \}$ be a Thurston geometry. 
Consider a translation triangle $A_1A_2A_3$ in the space $X$ whose translation surface defined in \cite{Cs-Sz25}, is denoted 
by $\cS_{A_1A_2A_3}^{X,t}$. Moreover, let the translation triangle $B_1B_2B_3$ be located in the same surface. We assume for the above two triangles that 
the translation curves 
$g_{A_1B_1}^X$, $g_{A_2B_2}^X$, $g_{C_1C_2}^X$ intersect at the proper point $S \subset \cS_{A_1A_2A_3}^{X,t}$. 
Then we say that the triangles $A_1A_2A_3$ and $B_1B_2B_3$ are {\it perspective with respect to the point, (perspective centrally)} in the 
$X$ sense. Moreover, if for the initial triangles it is true that ${M_{12}}=g_{A_1A_2}^X \cap g_{B_1B_2}^X$, $M_{13}=g_{A_1A_3}^X \cap 
g_{B_1B_3}^X$, $M_{23}=g_{A_2A_3}^X \cap g_{B_2B_3}^X$
proper points exist and fit on a translation curve, the triangles are {\it perspective with respect to a line, (perspective axially)} in the $X$.
\end{defn}
\begin{Theorem}[Desargues's theorem in $X$ space] 

If two triangles $A_1A_2A_3$ and $B_1B_2B_3$ are perspective centrally with respect to the point $S$, 
then the points $M_{12}$, $M_{13}$ and $M_{23}$ of intersection of the two triangles' corresponding
sides (translation curves) - if those intersection points exist - all lie on one translation curve. \qquad \qquad $\square$
\end{Theorem}
\begin{Theorem}[Pappus's hexagon theorem in $X$ space] 

Let $A_1, A_2, A_3 \subset \mathcal{S}^{X,t}_{g_a,g_b}$  
be three points on a translation curve $g_a^X$ and let $B_1, B_2, B_3 
\subset \mathcal{S}^{X,t}_{g_a,g_b}$ be three points on another translation 
curve $g_b^X$. 
If the translation curves $g_{A_1B_1}^X$, 
$g_{A_1B_1}^X$, $g_{A_2B_2}^X$, $g_{A_3B_3}^X$ intersect the translation 
curves $g_{A_3B_2}^X$,  $g_{A_1B_3}^X$, $g_{A_2B_1}^X$, respectively then 
then the three points of intersection $M_1=g_{A_1B_1}^X \cap g_{A_3B_2}^X$,
$M_2=g_{A_2B_2}^X \cap g_{A_1B_3}^X$, 
$M_3=g_{A_3B_3}^X \cap g_{A_2B_1}^X$ - if those 
intersection points exist - all lie on one translation curve. \qquad \qquad $\square$
\end{Theorem}
\begin{rem}
So the surface of the translation triangles of Thurston geometries defined in the work \cite{Cs-Sz25} is both Pappian and Desarguesian ``plane".
\end{rem}
Similar problems in homogeneous Thurston geometries
represent huge class of open mathematical problems. For
$\SXR$, $\HXR$, $\NIL$, $\SOL$, $\SLR$ geometries only very few results are known
(e.g. \cite{CsSz16}, \cite{Cs23}, \cite{MSzV},  \cite{MSz}, \cite{MSz12}, \cite{JJS}, \cite{Sz10}, \cite{Sz14-1}, \cite{Sz14-2}, \cite{Sz19}, \cite{Sz24}  ). 
Detailed studies are the objective of
ongoing research.

\vspace{0.5cm}

{\bf Funding} The authors declare that no funds, grants, or other support were
received during the preparation of this manuscript.
 
{\bf Data Availability Statement} The authors declare that this research does not
use data from any external sources or counterparts.

{\bf Declarations}

{\bf Competing Interest} The authors have no relevant financial or non-financial
interests to disclose.
%

%

\end{document}